\documentclass{amsart}
\usepackage{amssymb}
\usepackage{amsmath}
\usepackage{amsfonts}

\newtheorem{theorem}{Theorem}
\theoremstyle{plain}

\newtheorem{definition}{Definition}

\newtheorem{proposition}{Proposition}

\numberwithin{equation}{section}
\input{tcilatex}

\begin{document}
\title[On Continuous Fields of JB-algebras]{On Continuous Fields of
JB-algebras}
\author{Alexander A. Katz}
\address{Dr. Alexander A. Katz, Department of Mathematics and Computer
Science, St. John's College of Liberal Arts and Sciences, St. John's
University, 300 Howard Avenue, DaSilva Academic Center 314, Staten Island,
NY 10301, USA}
\email{katza@stjohns.edu}
\date{October 13-14, 2007}
\dedicatory{Dedicated to the memory of Professor George Bachman, Polytechnic
University, Brooklyn, NY, USA.}

\begin{abstract}
We introduce and study continuous fields of JB-algebras (which are real
non-associate analogues of C*-algebras). In particular, we show that for the
universal enveloping C*-algebra $C_{u}^{\ast }(B)$ for the JB-algebra $B$
defined by a continuous field of JB-algebras $A_{t},$ $t\in T,$\ on a
locally compact space $T$ there exists a decomposition of $C_{u}^{\ast }(B)$
into a continuous field of C*-algebras $C_{u}^{\ast }(A_{t}),$ $t\in T,$ on
the same space $T$, composed entirely of the universal enveloping
C*-algebras of the corresponding JB-algebras from the aforementioned
decomposition of the algebra $B$.
\end{abstract}

\maketitle

\section{Introduction and Preliminaries}

Banach associative regular *-algebras over $%
\mathbb{C}
$, so called \textit{C*-algebras}, were first introduces by Gelfand and
Naimark in the paper \cite{GelfandNaimark43}. Since then these algebras were
studied extensively by various authors. This theory now is a big subdomain
of the of Functional Analysis as a subject which found applications in
almost all branches of Modern Mathematics and Physics. For the basics of the
theory of C*-algebras, see for example Pedersen's monograph \cite{Pedersen79}%
. The basic theory of real associative analogues of C*-algebras, so called 
\textit{real C*-algebras}, is presented in Li's monograph \cite{Li03}.

In order to obtain a topological non-commutative version of Gelfand's
characterization of commutative C*-algebras, Dixmier and Douady in \cite%
{DixmierDouady63} introduced a notion of continuous fields of Banach spaces
and C*-algebras, which found important applications in classification of
C*-algebras (see \cite{Dixmier77}) and Theoretical Physics (see \cite%
{Landsman98}). According to them, \textbf{a continuous field of C*-algebras}%
\textit{\ }%
\begin{equation*}
(\mathfrak{B},\{\mathfrak{A}_{t},\varphi _{t}\}_{t\in T}),
\end{equation*}%
\textit{over a locally compact Hausdorff space }$T$\textit{\ consists of a
C*-algebra }$\mathfrak{B}$\textit{, a collection of C*-algebras }%
\begin{equation*}
\{\mathfrak{A}_{t}\}_{t\in T},
\end{equation*}%
\textit{and a set }%
\begin{equation*}
\{\varphi _{t}:\mathfrak{B}\rightarrow \mathfrak{A}_{t}\}_{t\in T},
\end{equation*}%
\textit{of surjective morphisms, such that:}

\textit{1). The function }%
\begin{equation*}
t\mapsto \left\Vert \varphi _{t}(x)\right\Vert ,
\end{equation*}%
\textit{is in }$C_{0}(T)$\textit{\ for all }$x\in \mathfrak{B};$

\textit{2). The norm of any }$x\in \mathfrak{B}$\textit{\ is }%
\begin{equation*}
\left\Vert x\right\Vert _{\mathfrak{B}}=\underset{t\in T}{\sup }\left\Vert
\varphi _{t}(x)\right\Vert _{\mathfrak{A}_{t}};
\end{equation*}

\textit{3). For any }$f\in C_{0}(T)$\textit{\ and }$x\in \mathfrak{B},$%
\textit{\ there is an element }%
\begin{equation*}
fx\in \mathfrak{B},
\end{equation*}%
\textit{for which }%
\begin{equation*}
\varphi _{t}(fx)=f(t)\varphi _{t}(x),
\end{equation*}%
\textit{for all }$t\in T.$

\textit{A \textbf{section} of the field is an element }$\{x_{t}\}_{t\in T}$%
\textit{\ of }%
\begin{equation*}
\dprod\limits_{t\in T}\mathfrak{A}_{t},
\end{equation*}%
\textit{for which there is an element }$x\in \mathfrak{B}$\textit{\ such
that }%
\begin{equation*}
x_{t}=\varphi _{t}(x),
\end{equation*}%
\textit{for all }$t\in T.$ 

\textit{One can see that }$\mathfrak{B}$\textit{\ can be identified with the
space of sections of the field, seen as a C*-algebra under pointwise scalar
multiplication, addition, adjointing, and operator multiplication, by means }%
\begin{equation*}
\{\varphi _{t}(x)\}_{t\in T}\leftrightarrow x.
\end{equation*}%
\textit{In particular, }%
\begin{equation*}
x=y,
\end{equation*}%
\textit{iff }%
\begin{equation*}
\varphi _{t}(x)=\varphi _{t}(y),
\end{equation*}%
\textit{for all }$t.$\textit{\ It is natural that algebra }$\mathfrak{B}$%
\textit{\ is called a the \textbf{C*-algebra of the continuous field of
C*-algebras}.}

Since the beginning of the theory of complex C*-algebras, there were
numerous attempts to extend this theory to non-associative algebras which
are close to associative, in particular to Jordan algebras. In fact, Alfsen,
Shultz and St\o rmer in \cite{AlfsenShultzStoermer78} defined so called 
\textit{JB-algebras} as the real Banach--Jordan algebras satisfying for all
pairs of elements $x$ and $y$ the inequality of fineness 
\begin{equation*}
\left\Vert x^{2}+y^{2}\right\Vert \geq \left\Vert x\right\Vert ^{2},
\end{equation*}%
and regularity condition 
\begin{equation*}
\left\Vert x^{2}\right\Vert =\left\Vert x\right\Vert ^{2}.
\end{equation*}%
If $\mathfrak{A}$ is a C*--algebra, or a real C*-algebra, then the
self-adjoint part $\mathfrak{A}_{sa}$ of $\mathfrak{A}$ is a JB-algebra
under the Jordan product 
\begin{equation*}
x\circ y=\frac{(xy+yx)}{2}.
\end{equation*}%
Closed subalgebras of $\mathfrak{A}_{sa}$, for some C*-algebra or real
C*-algebra $\mathfrak{A}$, become relevant examples of JB-algebras, and are
called \textit{JC-algebras}.

The basic theory of JB-algebras is fully treated in monograph of
Hanche-Olsen and St\o rmer \cite{Hanche-OlsenStoermer84}.In particular, in
this monograph there is the following theorem which was for the first time
presented by Alfsen, Hanche-Olsen and Shultz in the paper \cite%
{AlfsenHanche-OlsenShultz80}.

\begin{theorem}[Alfsen, Hanche-Olsen, Shultz \protect\cite%
{AlfsenHanche-OlsenShultz80}]
For an arbitrary JB-algebra $A$ there exists a unique up to an isometric
*-isomorphism a C*-algebra $C_{u}^{\ast }(A)$ (\textbf{the universal
enveloping C*-algebra for the JB-algebra }$\mathbf{A}$), and a Jordan
homomorphism 
\begin{equation*}
\psi _{A}:A\rightarrow C_{u}^{\ast }(A)_{sa},
\end{equation*}%
from $A$ to the self-adjoint part of $C_{u}^{\ast }(A),$ such that:

1). $\psi _{A}(A)$ generates $C_{u}^{\ast }(A)$ as a C*-algebra;

2). for any pair composed of a C*-algebra $\mathfrak{A}$ and a Jordan
homomorphism 
\begin{equation*}
\rho :A\rightarrow \mathfrak{A}_{sa},
\end{equation*}%
from A into the self-adjoint part of $\mathfrak{A}$, there exists a
*-homomorphism 
\begin{equation*}
\widehat{\rho }:C_{u}^{\ast }(A)\rightarrow \mathfrak{A,}
\end{equation*}%
from the C*-algebra $C_{u}^{\ast }(A)$ into C*-algebra $\mathfrak{A}$, such
that 
\begin{equation*}
\rho =\widehat{\rho }\circ \psi _{A};
\end{equation*}

3). there exists a *-antiautomorphism $\Phi $ of order 2 on the C*-algebra $%
C_{u}^{\ast }(A)$, such that 
\begin{equation*}
\Phi (\psi _{A}(a))=\psi _{A}(a),
\end{equation*}%
$\forall a\in A.$ \ \ \ \ \ \ \ \ \ \ \ \ \ \ \ \ \ \ \ \ \ \ \ \ \ \ \ \ \
\ \ \ \ \ \ \ \ \ \ \ \ \ \ \ \ \ \ \ \ \ \ \ \ \ \ \ \ \ \ \ \ \ \ \ \ \ \
\ \ \ \ \ \ \ \ \ \ \ \ \ \ \ \ \ \ \ \ \ $\square $
\end{theorem}

Our plan is to define a \textit{continuous field of JB-algebras}, the 
\textit{JB-algebra of the continuous field of JB-algebras}, and be able in
the spirit of Theorem 1 above to associate in a universal sense with each
JB-algebra of the continuous field of JB-algebras a C*-algebra of the
continuous field of C*-algebras.

\section{Continuous fields of JB-algebras}

Let us first introduce a continuous field of JB-algebras.

\begin{definition}
A\textbf{\ continuous field of JB-algebras }%
\begin{equation*}
(B,\{A_{t},\varphi _{t}\}_{t\in T}),
\end{equation*}%
\textit{over a locally compact Hausdorff space }$T$\textit{\ consists of a
JB-algebra }$B$\textit{, a collection of JB-algebras }$\{A_{t}\}_{t\in T},$%
\textit{and a set }%
\begin{equation*}
\{\varphi _{t}:B\rightarrow A_{t}\}_{t\in T},
\end{equation*}%
\textit{of surjective morphisms, such that: }

\textit{1). The function }%
\begin{equation*}
t\mapsto \left\Vert \varphi _{t}(x)\right\Vert ,
\end{equation*}%
\textit{is in }$C_{0}(T)$\textit{\ for all }$x\in B;$ 

\textit{2). The norm of any }$x\in B$\textit{\ is }%
\begin{equation*}
\left\Vert x\right\Vert _{B}=\underset{t\in T}{\sup }\left\Vert \varphi
_{t}(x)\right\Vert _{A_{t}};
\end{equation*}

\textit{3). For any }$f\in C_{0}(T)$\textit{\ and }$x\in B,$\textit{\ there
is an element }%
\begin{equation*}
fx\in B,
\end{equation*}%
\textit{for which }%
\begin{equation*}
\varphi _{t}(fx)=f(t)\varphi _{t}(x),
\end{equation*}%
\textit{for all }$t\in T.$
\end{definition}

\textit{A \textbf{section} of the field is an element }$\{x_{t}\}_{t\in T}$%
\textit{\ of }%
\begin{equation*}
\dprod\limits_{t\in T}A_{t},
\end{equation*}%
\textit{for which there is an element }$x\in B$\textit{\ such that }%
\begin{equation*}
x_{t}=\varphi _{t}(x),
\end{equation*}%
\textit{for all }$t\in T.$ 

\textit{We identify }$B$\textit{\ with the space of sections of the field,
seen as a JB-algebra under pointwise scalar multiplication, addition,
operator multiplication, by means }%
\begin{equation*}
\{\varphi _{t}(x)\}_{t\in T}\leftrightarrow x.
\end{equation*}%
\textit{In particular, }%
\begin{equation*}
x=y,
\end{equation*}%
\textit{iff }%
\begin{equation*}
\varphi _{t}(x)=\varphi _{t}(y),
\end{equation*}%
\textit{for all }$t.$\textit{\ It is natural that algebra }$B$\textit{\ is
called a the \textbf{JB-algebra of the continuous field of JB-algebras}.}

Now we will establish a few properties of the continuous field of
JB-algebras. The first one is about locally uniform closedness of the
continuous field of JB-algebras.

\begin{proposition}
The JB-algebra B of sections of a continuous field of JB-algebras is \textbf{%
locally uniformly closed}, i.e. if 
\begin{equation*}
x\in \dprod\limits_{t\in T}A_{t},
\end{equation*}%
is such that for every 
\begin{equation*}
s\in T,
\end{equation*}%
and every $\varepsilon >0$ there exists 
\begin{equation*}
y_{s}\in B,
\end{equation*}%
and a neighborhood 
\begin{equation*}
V_{s}\subset T,
\end{equation*}%
of $s$ in which 
\begin{equation*}
\left\Vert x_{t}-\varphi _{t}(y_{s})\right\Vert <\varepsilon ,
\end{equation*}%
$\ $for all 
\begin{equation*}
t\in V_{s},
\end{equation*}%
and also 
\begin{equation*}
\underset{t\rightarrow \infty }{\lim }\left\Vert x_{t}\right\Vert =0,
\end{equation*}%
then%
\begin{equation*}
x\in B.
\end{equation*}

Alternatively, if the function 
\begin{equation*}
t\mapsto \left\Vert x_{t}-z_{t}\right\Vert ,
\end{equation*}%
lies in $C_{0}(T)$ for each 
\begin{equation*}
z\in B,
\end{equation*}%
then 
\begin{equation*}
x\in B.
\end{equation*}
\end{proposition}

\begin{proof}
Under conditions of the first part of the Proposition, there exists a
compact set 
\begin{equation*}
K\subseteq T,
\end{equation*}%
for which 
\begin{equation*}
\left\Vert x_{t}\right\Vert <\varepsilon ,
\end{equation*}%
outside of $K,$ as well as a finite cover 
\begin{equation*}
\{V_{t_{1}},...,V_{t_{n}}\},
\end{equation*}%
of $K$, 
\begin{equation*}
K\subseteq \{V_{t_{1}},...,V_{t_{n}}\}.
\end{equation*}

Now we have to recall a notion of a partition of unity on $K$ subordinate to
this cover (see for example \cite{Pedersen89} and \cite{Landsman98}). Let $K$
be a Hausdorff space, and 
\begin{equation*}
\{V_{\alpha }\}_{\alpha \in \Lambda },
\end{equation*}%
be a \textbf{locally finite open cover} of $K,$ i.e. each point of of $K$
has a neighborhood that intersects only with a finite number of the sets $%
V_{\alpha }.$ A \textbf{partition of unity subordinate to the given cover}
is a collection of positive functions 
\begin{equation*}
\{u_{\alpha }\}_{\alpha \in \Lambda },
\end{equation*}%
such that each $u_{\alpha }$ is a compactly supported continuous real-valued
function with 
\begin{equation*}
\dsum\limits_{\alpha \in \Lambda }u_{\alpha }=1.
\end{equation*}%
A partition of unity always exists when $K$ is paracompact (see \cite%
{Pedersen89}). So, let us take a partition of unity 
\begin{equation*}
\{u_{i}\}_{i=1}^{n},
\end{equation*}%
on $K$ subordinate to the aforementioned finite cover 
\begin{equation*}
\{V_{t_{1}},...,V_{t_{n}}\}.
\end{equation*}%
Let us consider the 
\begin{equation*}
y=\dsum\limits_{i=1}^{n}u_{i}y_{t_{i}}.
\end{equation*}%
From Definition 1.3 it follows that 
\begin{equation*}
y\in B,
\end{equation*}%
and satisfies the condition 
\begin{equation*}
\underset{t\in T}{\sup }\left\Vert x_{t}-y_{t}\right\Vert <\varepsilon .
\end{equation*}%
Therefore, from Definition 1.2 and completeness of $B$ it follows that 
\begin{equation*}
x\in B.
\end{equation*}

Now, given any 
\begin{equation*}
x\in \dprod\limits_{t\in T}A_{t},
\end{equation*}%
and 
\begin{equation*}
s\in T,
\end{equation*}%
because $\varphi _{s}$ is surjective, there exists an element 
\begin{equation*}
y_{s}\in B,
\end{equation*}%
such that 
\begin{equation*}
x_{s}=\varphi _{s}(y_{s}).
\end{equation*}%
The assumption of the second part of Proposition 1 then implies that the
conditions in the first part of this Proposition are satisfied, such that 
\begin{equation*}
x\in B.
\end{equation*}
\end{proof}

The following Proposition gives conditions for the existence and uniqueness
of a continuous field of JB-algebras whose collection of sections contains a
subset possessing some natural properties.

\begin{proposition}
Let 
\begin{equation*}
\{A_{t}\}_{t\in T},
\end{equation*}%
be a family of JB-algebras indexed by a locally compact Hausdorff space $T$,
and a subset 
\begin{equation*}
\widetilde{B}\subseteq \dprod\limits_{t\in T}A_{t},
\end{equation*}%
that satisfies the following properties:

1). The set 
\begin{equation*}
\{x_{t}:x\in \widetilde{B}\},
\end{equation*}%
is dense in $A_{t}$ for each $t\in T;$

2). The function 
\begin{equation*}
t\mapsto \left\Vert x_{t}\right\Vert ,
\end{equation*}%
lies in $C_{0}(T)$ for each $x\in \widetilde{B};$

3). The set $\widetilde{B}$ is a Jordan algebra under pointwise operations.

Then there exists a unique continuous field of JB-algebras 
\begin{equation*}
(B,\{A_{t},\varphi _{t}\}_{t\in T}),
\end{equation*}%
whose collection of sections contains $\widetilde{B}.$ Namely, $B$ consists
of all 
\begin{equation*}
x\in \dprod\limits_{t\in T}A_{t},
\end{equation*}%
for which the function 
\begin{equation*}
x\mapsto \left\Vert x_{t}-z_{t}\right\Vert ,
\end{equation*}%
lies in $C_{0}(T)$ for each $z\in \widetilde{B},$ regarded as JB-algebra
under pointwise operations, and the norm of Definition 1.2. Finally, 
\begin{equation*}
\varphi _{t}(x)=x_{t},
\end{equation*}%
$t\in T,$  is the evaluation map.
\end{proposition}

\begin{proof}
We show first that the algebra $B$ defined above is locally uniformly
closed. With the objects $x,s,\varepsilon ,y_{s}$ and V as specified in
Proposition 1, take $z\in \widetilde{B}$ arbitrary, and define the functions 
\begin{equation*}
f_{xz}:t\mapsto \left\Vert x_{t}-z_{t}\right\Vert ,
\end{equation*}%
and 
\begin{equation*}
f_{yz}:t\mapsto \left\Vert \varphi _{t}(y_{s})-z_{t}\right\Vert .
\end{equation*}%
Using the triangle inequality for the norm in Banach space, we get that 
\begin{equation*}
\left\vert (\left\Vert x\right\Vert -\left\Vert y\right\Vert )\right\vert
\leq \left\Vert x-y\right\Vert ,
\end{equation*}%
and that gives us that 
\begin{equation*}
\left\vert f_{xz}(t)-f_{yz}(t)\right\vert <\varepsilon ,
\end{equation*}%
for all $t\in V.$ By assumption, the function $f_{yz}$ is continuous, so
that 
\begin{equation*}
\left\vert f_{yz}(t)-f_{yz}(s)\right\vert <\varepsilon ,
\end{equation*}%
for all $t$'s in some neighborhood $V^{\prime }$ of $s.$ Combining the
inequalities, we get 
\begin{equation*}
\left\vert f_{xz}(t)-f_{xz}(s)\right\vert <3\varepsilon ,
\end{equation*}%
for all 
\begin{equation*}
t\in V\cap V^{\prime }.
\end{equation*}%
Therefore $f_{xz}$ is continuous at $s,$ which was arbitrary, so that $x\in B
$ by the definition of $B$.

Now, we show uniqueness of $B$. Using this property one can easily see that $%
B$ is a JB-algebra, and that the condition 3 in Definition 1 is satisfied.
It is clear from Definition 1.1 and the definition of $B$ in Proposition 2
that $B$ is maximal. On the other hand, according to the second part of
Proposition 1, $B$ is minimal, so, $B$ is unique.
\end{proof}

We are ready now to present the main result of the paper.

\begin{definition}
An *-isomorphism (resp. Jordan isomorphism) of continuous fields of
C*-algebras (resp. JB-algebras) over the same base Hausdorff locally compact
topological space T is the isometric *-isomorphism (Jordan isometric
isomorphism) of the C*-algebras (resp. JB-algebras) of the continuous fields
via a map respecting the fibers.
\end{definition}

\begin{theorem}
For an arbitrary continuous field of JB-algebras\textbf{\ }%
\begin{equation*}
(B,\{A_{t},\varphi _{t}\}_{t\in T}),
\end{equation*}%
\textit{over a locally compact Hausdorff space }$T,$ there exists a unique
up to an *-isomorphism a continuous field of C*-algebras 
\begin{equation*}
(C_{u}^{\ast }(B),\{C_{u}^{\ast }(A_{t}),\widehat{\varphi }_{t}\}_{t\in T}),
\end{equation*}%
(\textbf{the universal enveloping continuous field of C*-algebras for the
continuous field of JB-algebras }%
\begin{equation*}
(B,\{A_{t},\varphi _{t}\}_{t\in T})
\end{equation*}%
\textbf{over the same base space }$T$), and a Jordan homomorphism 
\begin{equation*}
\psi _{B}:B\rightarrow C_{u}^{\ast }(B)_{sa},
\end{equation*}%
from $B$ to the self-adjoint part of $C_{u}^{\ast }(B),$ as well as a family
of Jordan homomorphisms 
\begin{equation*}
\psi _{A_{t}}:A_{t}\rightarrow C_{u}^{\ast }(A_{t})_{sa},
\end{equation*}%
$t\in T,$ from $A_{t}$ to the self-adjoint part of $C_{u}^{\ast }(A_{t}),$
for each $t\in T,$ such that:

1). $\psi _{B}(B)$ generates $C_{u}^{\ast }(B)$ as a C*-algebra, and each $%
\psi _{A_{t}}(B)$ generates each $C_{u}^{\ast }(A_{t})$ as a C*-algebra for
each $t\in T$;

2). for any pair composed of a continuous field of C*-algebras $(\mathfrak{B}%
,\{\mathfrak{A}_{t},\widetilde{\varphi }_{t}\}_{t\in T})$ and a family of
Jordan homomorphisms 
\begin{equation*}
\rho :B\rightarrow \mathfrak{B}_{sa},
\end{equation*}%
from $B$ into the self-adjoint part of $\mathfrak{B}$, and 
\begin{equation*}
\rho _{t}:A_{t}\rightarrow (\mathfrak{A}_{t})_{sa},
\end{equation*}%
for each $t\in T,$\ from $A_{t}$ into the self-adjoint part of $\mathfrak{A}%
_{t},$ there exist a *-homomorphism 
\begin{equation*}
\widehat{\rho }:C_{u}^{\ast }(B)\rightarrow \mathfrak{B},
\end{equation*}%
from the C*-algebra $C_{u}^{\ast }(B)$ into C*-algebra $\mathfrak{B}$, and a
family of *-homomorphisms 
\begin{equation*}
\widehat{\rho }_{t}:C_{u}^{\ast }(A_{t})\rightarrow \mathfrak{A}_{t},
\end{equation*}%
for each $t\in T,$ from the C*-algebra $C_{u}^{\ast }(A_{t})$ into
C*-algebra $\mathfrak{A}_{t}$ such that 
\begin{equation*}
\rho =\widehat{\rho }\circ \psi _{B},
\end{equation*}%
and 
\begin{equation*}
\rho _{t}=\widehat{\rho }_{t}\circ \psi _{A_{t}},
\end{equation*}%
for each $t\in T;$

3). there exists a *-antiautomorphism $\Phi $ of order 2 on the C*-algebra $%
C_{u}^{\ast }(B)$, such that 
\begin{equation*}
\Phi (\psi _{B}(x))=\psi _{B}(x),
\end{equation*}%
$\forall x\in B,$ as well as there exists a family of  *-antiautomorphism $%
\Phi _{t}$ of order 2 on the C*-algebra $C_{u}^{\ast }(A_{t})$, for each $%
t\in T,$ such that 
\begin{equation*}
\Phi (\psi _{A_{t}}(x_{t}))=\psi _{A_{t}}(x_{t}),
\end{equation*}%
$\forall x_{t}\in A_{t},$ and every $t\in T.$
\end{theorem}

\begin{proof}
Let $(B,\{A_{t},\varphi _{t}\}_{t\in T})$ be a given continuous field of
JB-algebras. Let $C_{u}^{\ast }(B)$ be the universal enveloping C*-algebra
for the JB-algebra of the continuous field, and the family of C*-algebras $%
C_{u}^{\ast }(A_{t})$ for each $t\in T$ be the universal enveloping
C*-algebra for the JB-algebra $A_{t},$ $t\in T.$ Let 
\begin{equation*}
\psi _{B}:B\rightarrow C_{u}^{\ast }(B)_{sa},
\end{equation*}%
and for each $t\in T,$ 
\begin{equation*}
\varphi _{t}:B\rightarrow A_{t},
\end{equation*}%
and 
\begin{equation*}
\psi _{A_{t}}:A_{t}\rightarrow C_{u}^{\ast }(A_{t})_{sa}.
\end{equation*}%
From Theorem 1 it follows that $\psi _{B}(B)$ is dense in $C_{u}^{\ast
}(B)_{sa},$ and $\psi _{A_{t}}(A_{t})$ is dense in $C_{u}^{\ast }(A_{t})_{sa}
$ for each $t\in T$. So, without a loss of generality we can assume that for
each 
\begin{equation*}
x,y\in C_{u}^{\ast }(B)_{sa}
\end{equation*}%
there exist $a_{n},b_{n}\in B,$ such that 
\begin{equation*}
x=\underset{n\rightarrow \infty }{\lim }\psi _{B}(a_{n}),
\end{equation*}%
and 
\begin{equation*}
y=\underset{n\rightarrow \infty }{\lim }\psi _{B}(b_{n}),
\end{equation*}%
where the limit is taken in the norm of $C_{u}^{\ast }(B),$ as well as 
\begin{equation*}
x_{t}=\widehat{\varphi }_{t}(x)=\underset{n\rightarrow \infty }{\lim }\psi
_{A_{t}}(\varphi _{t}(a_{n})),
\end{equation*}%
and 
\begin{equation*}
y_{t}=\widehat{\varphi }_{t}(y)=\underset{n\rightarrow \infty }{\lim }\psi
_{A_{t}}(\varphi _{t}(b_{n})),
\end{equation*}%
$t\in T$, where the limit is taken in the norm of $C_{u}^{\ast }(A_{t}).$
For each $t\in T,$ we will define 
\begin{equation*}
\widehat{\varphi }_{t}:C_{u}^{\ast }(B)\rightarrow C_{u}^{\ast }(A_{t}),
\end{equation*}%
the following way: 
\begin{equation*}
\widehat{\varphi }_{t}(x+iy)=x_{t}+iy_{t}.
\end{equation*}%
Because $C_{u}^{\ast }(A_{t})_{sa}$ is norm closed for each $t\in T$ (see 
\cite{Pedersen79}), the last identity it well defined. Moreover, from the
fact that $\varphi _{t}$ is surjective for each $t\in T$ it follows that $%
\widehat{\varphi }_{t}$ is surjective for each $t\in T$ as well. Thus, 
\begin{equation*}
(C_{u}^{\ast }(B),\{C_{u}^{\ast }(A_{t}),\widehat{\varphi }_{t}\}_{t\in T}),
\end{equation*}%
is in fact a continuous field of C*-algebras. The rest of the Theorem is
obtained as a corollary by application of Theorem 1 in fibers, and
Propositions 1 and 2.
\end{proof}


\begin{thebibliography}{99}
\bibitem{AlfsenHanche-OlsenShultz80} \textbf{Alfsen, E.M.; Hanche-Olsen, H.;
Shultz, F.W.,} \textit{State spaces of C*-algebras.} (English) Acta Math.
Vol. 144 (1980), No. 3-4, pp. 267-305.

\bibitem{AlfsenShultzStoermer78} \textbf{Alfsen, E.M.; Shultz, F.W.; St\o %
rmer, E.,} \textit{A Gelfand-Naimark theorem for Jordan algebras.} (English)
Advances in Math. Vol. 28 (1978), No. 1, pp. 11-56.

\bibitem{Dixmier77} \textbf{Dixmier, J.,} \textit{C*-algebras.} (English,
Translated from French) North-Holland Mathematical Library, Vol. 15.
North-Holland Publishing Co., Amsterdam-New York-Oxford, 492 pp., (1977).

\bibitem{DixmierDouady63} \textbf{Dixmier, J.; Douady, A.,} \textit{Champs
Continus d'Espaces Hilbertiens et de C*-Algebres.} (French) \textit{Bull.
Soc. math. France} Vol. 91 (1963), pp. 227--284.

\bibitem{GelfandNaimark43} \textbf{Gelfand, I.M.; Naimark, M.A.,} \textit{On
the embedding of normed rings into the ring of operators in Hilbert space.}
(English. Russian summary) Rec. Math. [Mat. Sbornik] N.S. Vol. 12(54)
(1943), pp. 197-213.

\bibitem{Hanche-OlsenStoermer84} \textbf{Hanche-Olsen, H.; St\o rmer, E.,} 
\textit{Jordan operator algebras.} (English), Monographs and Studies in
Mathematics, Vol. 21. Boston - London - Melbourne: Pitman Advanced
Publishing Program. VIII, 183 pp., (1984).

\bibitem{Landsman98} \textbf{Landsman, N.P.,} \textit{Mathematical topics
between classical and quantum mechanics.} (English), Springer Monographs in
Mathematics. Springer-Verlag, New York, 529 pp., (1998).

\bibitem{Li03} \textbf{Li, B.,} \textit{Real operator algebras.} (English)
World Scientific Publishing Co., Inc., River Edge, NJ, 241 pp., (2003).

\bibitem{Pedersen79} \textbf{Pedersen, G.K.,} \textit{C*-algebras and their
automorphism groups.} (English), London Mathematical Society Monographs.
Vol. 14. London - New York -San Francisco: Academic Press., 416 pp., (1979).

\bibitem{Pedersen89} \textbf{Pedersen, G.K.,} \textit{Analysis now.}
(English) Graduate Texts in Mathematics, Vol. 118. Springer-Verlag, New
York, 277 pp., (1989).
\end{thebibliography}
\end{document}